\documentclass[12pt]{article}
\usepackage{latexsym}
\usepackage{amsfonts}
\usepackage{amsmath} 
\usepackage{geometry}
\usepackage{amssymb}
\usepackage[mathscr]{euscript}
\usepackage[cal=boondoxo]{mathalfa} % mathcal
\usepackage{enumerate}
\makeatletter \@addtoreset{equation}{section} \makeatother \newtheorem{lemma}[equation]{{\sc Lemma}}\newtheorem{corollary}[equation]{{\sc Corollary}}\newtheorem{proposition}[equation]{{\sc Proposition}}\newtheorem{exampleT}[equation]{{\sc Example}}\newtheorem{exerciseT}[equation]{{\sc Exercise}}\newtheorem{remarkT}[equation]{{\sc Remark}}\newtheorem{definitionT}[equation]{{\sc Definition}}\newtheorem{theorem}[equation]{{\sc Theorem}}\renewcommand{\theequation}{\thesection.\arabic{equation}}

\def \qed {{\nopagebreak}\hfill \vrule height3pt width 3pt depth 0pt \vskip 3mm}

\newcommand{\definition}[1]{\begin{definitionT}{\sf  #1}\end{definitionT}}

\newcommand{\cld}{{\mathcal D}}

\newcommand{\clf}{{\mathcal F}}

\newcommand{\clp}{{\mathcal P}}

\newcommand{\Eb}{{{\sf E}}}
\newcommand{\Pb}{{{\sf P}}}

\newcommand{\Qb}{{{\sf Q}}}

\newcommand{\C}{{\mathbb C}}

\newcommand{\E}{{\mathbb E}}
\newcommand{\F}{{\mathbb F}}
\newcommand{\G}{{\mathbb G}}

\newcommand{\K}{{\mathbb K}}
\renewcommand{\L}{{\mathbb L}}
\newcommand{\M}{{\mathbb M}}

\def\Ind {\text{\large 1}}

\newcommand{\pr}{\NI{\sc Proof : }}
\newcommand{\NI}{\noindent}
\newcommand{\bea}{\begin{eqnarray}}
\newcommand{\eea}{\end{eqnarray}}
\newcommand{\be}{\begin{equation}}
\newcommand{\ee}{\end{equation}}
\newcommand{\ben}{\begin{eqnarray*}}
\newcommand{\een}{\end{eqnarray*}}
\newcommand{\bt}{\begin{theorem}}
\newcommand{\et}{\end{theorem}}
\newcommand{\bl}{\begin{lemma}}
\newcommand{\el}{\end{lemma}}
\newcommand{\ber}{\begin{exercise}}
\newcommand{\eer}{\end{exercise}}
\newcommand{\bem}{\begin{example}}
\newcommand{\eem}{\end{example}}
\newcommand{\bc}{\begin{corollary}}
\newcommand{\ec}{\end{corollary}}
\newcommand{\bp}{\begin{proposition}}
\newcommand{\ep}{\end{proposition}}

\renewcommand{\d}{{\hspace{0.6pt}d\hspace{0.1pt}}}

\def\title #1{\begin{center}
{\large {\bf #1}}
\end{center}}
\def\place #1{\centerline {\it #1}}

\begin{document}

\begin{titlepage}
\vspace*{4mm}
\title{On the Second Fundamental Theorem of Asset Pricing}
\vskip 2em
\centerline{\mbox{{\sc Rajeeva L. Karandikar}\hspace{10pt} and
\hspace{10pt}{\sc B. V. Rao
 }}}
\vskip 1em
\place{Chennai Mathematical Institute, Chennai.}
\vskip 9em
\centerline{{\sc Abstract}}

{ Let $X^1,\ldots, X^d$ be  sigma-martingales on $(\Omega,\clf, \Pb)$. We show that every bounded martingale (with respect to the underlying filtration) admits an integral representation w.r.t. $X^1,\ldots, X^d$ if and only if there is no equivalent probability measure (other than $\Pb$) under which $X^1,\ldots,X^d$ are sigma-martingales.

From this we deduce the second fundamental theorem of asset pricing-  that completeness of a market is equivalent to uniqueness of Equivalent Sigma-Martingale Measure (ESMM).}
\vfill
\hrule width100pt height .8pt

\vskip 2mm
\noindent {\it 2010 Mathematics Subject Classication:}  91G20, 60G44, 97M30, 62P05.\\
\noindent {\it Key words and phrases:} Martingales, Sigma Martingales, Stochastic Calculus, Martingale Representation, No Arbitrage, Completeness of Markets.
\end{titlepage}

\section{Introduction}
The (first) fundamental theorem of asset pricing says that a market consisting of finitely many stocks satisfies the {\em No Arbitrage property} (NA) if and only there exists an {\em Equivalent Martingale Measure} (EMM)- {\em i.e. there exists an equivalent probability measure under which the (discounted) stocks are (local) martingales.} The No Arbitrage property has to be suitably defined when we are dealing in continuous time, where one rules out approximate arbitrage in the class of admissible strategies. For a precise statement in the case when the underlying processes are locally bounded, see Delbaen and Schachermayer \cite{DS}. Also see Bhatt and Karandikar \cite{BK15} for an alternate formulation, where the approximate arbitrage is defined only in terms of simple strategies. For the general case, the result is true when local martingale in the statement above is replaced by sigma-martingale. See Delbaen and Schachermayer \cite{DS98}. They have an example where the No Arbitrage property holds but there is no equivalent measure under which the underlying process is a local martingale. However, there is an equivalent measure under which the process is a sigma-martingale. 
 
The second fundamental theorem of asset pricing says that the market is complete ({\em i.e. every contingent claim can be replicated by trading on the underlying securities}) if and only if the EMM is unique. Interestingly, this property was studied by probabilists well before the connection between finance and stochastic calculus was established (by Harrison-Pliska \cite{HP}). The  completeness of market is same as the question: when is every martingale representable as a stochastic integral w.r.t. a given set of martingales $\{M^1,\ldots ,M^d\}$. When  $M^1,\ldots ,M^d$ is the $d$-dimensional Wiener Process, this property was proven by Ito \cite{Ito}. Jacod and Yor \cite{JY} proved that if $M$ is a $\Pb$-local martingale, then every martingale $N$ admits a representation as a stochastic integral w.r.t. $M$ if and only if there is no probability measure $\Qb$ (other than $\Pb$) such that  $\Qb$ is equivalent to $\Pb$ and $M$ is a $\Qb$-local martingale. The situation in higher dimension is more complex. The obvious generalisation to higher dimension is not true as was noted by Jacod-Yor \cite{JY}.

To remedy the situation, a notion of vector stochastic integral was introduced- where a vector valued predictable process is the integrand and vector valued martingale is the integrator. The resulting integrals yield a class larger than the linear space generated by component wise integrals. See \cite{J80}, \cite{CS}. However, one has to prove various properties of stochastic integrals once again. 

Here we achieve the same objective in another fashion avoiding defining integration again from scratch. In the same breath, we also take into account the general case, when the underlying processes need not be bounded but satisfy the property NFLVR and thus one has a equivalent sigma-martingale measure (ESMM). We say that a martingale $M$ admits a integral representation w.r.t. $(X^1,X^2,\ldots, X^d)$ if there exits predictable $f, g^j$ such that $g^j\in\L(X^j)$,
\[Y_t=\sum_{j=1}^d\int_0^t g^j_s\d X^j_s,\]
$f\in\L(Y)$ and
\[M_t=\int_0^t f_s\d Y_s.\]
The security $Y$ can be thought of as a mutual fund or an index fund and the investor is trading on such a fund trying to replicate the security $M$. 

We will show that if for a multidimensional r.c.l.l. process $(X^1,X^2,\ldots, X^d)$ an ESMM exists, then all bounded martingales admit a representation w.r.t $X^j$, $1\le j\le d$ if and only if ESMM is unique.

\section{Preliminaries and Notation}
Let us start with some notations. $ (\Omega, \clf, \Pb)$ denotes a complete probability space with a filtration $ (\clf_\centerdot)=  \{\clf_t:\;t\ge 0 \}$ such that $\clf_0$ consists of all $\Pb$-null sets (in $\clf$) and  
\[\cap_{t>s}\clf_t=\clf_s \;\;\forall s\ge 0.\]

For various notions, definitions and standard results on stochastic integrals, we refer the reader to Jacod \cite{J78} or Protter \cite{P}.

For a local martingale $M$, let $\L^1_m(M)$ be the class of predictable processes $f$ such that there exists a sequence of stopping times $\sigma_k\uparrow\infty$ with
\[\Eb[\{\int_0^{\sigma_k}f^2_s\d [M,M]_s\}^\frac{1}{2}]<\infty\]
and for such an $f$, $N=\int f\d M$ is defined and is a local martingale. 

Let $\M$ denote the class of martingales and  for $M^1,M^2,\ldots , M^d\in\M$  let
\[\C(M^1,M^2,\ldots , M^d)=\{Z\in\M\,:\,Z_t=Z_0+\sum_{j=1}^d\int_0^t f^j_s\d M^j_s,\;f^j\in\L^1_m(M^j)\}\]
and for $T<\infty$, let
\[\tilde{\K}_T(M^1,M^2,\ldots , M^d)=\{Z_T\,:\, Z\in \C(M^1,M^2,\ldots , M^d)\}.\]

For the case $d=1$, Yor \cite{Yor} had proved that $\tilde{\K}_T$ is a closed subspace of $\L^1(\Omega,\clf,\Pb)$. The problem in case $d>1$ is that in general $\tilde{\K}_T(M^1,M^2,\ldots , M^d)$ need not be closed. Jacod-Yor \cite{JY} gave an example where $M^1,M^2$ are continuous square integrable martingales and $\tilde{\K}_T(M^1,M^2)$ is not closed.

For martingales $M^1,M^2,\ldots ,M^d$, let
\[\F(M^1,\ldots , M^d)=\{Z\in\M\,: Z_t=Z_0+\int_0^t f_s\d Y_s, \,f\in\L^1_m(Y),\,Y\in  \C(M^1,\ldots , M^d)\}.\]
Let
\[\K_T(M^1,M^2,\ldots , M^d)=\{Z_T\,:\,Z\in\F(M^1,M^2,\ldots , M^d)\}.\]
The main result of the next section is
\bt \label{aztm1} Let $M^1,M^2,\ldots ,M^d$ be martingales. Then $\K_T(M^1,M^2,\ldots , M^d)$ is closed in $\L^1(\Omega,\clf, \Pb)$. \et

This would be deduced from
\bt \label{aztm2}  Let $M^1,M^2,\ldots ,M^d$ be  martingales and $Z^n\in \F(M^1,M^2,\ldots , M^d)$ be such that $\Eb[\lvert Z^n_t-Z_t\rvert ]\rightarrow 0$ for all $t$. Then $Z\in \F(M^1,M^2,\ldots , M^d)$.
\et

When $M^1,M^2,\ldots ,M^d$ are square integrable martingales, the analogue of Theorem \ref{aztm1} for $\L^2$ follows from the work of Davis-Varaiya \cite{DV}. However, for the EMM characterisation via integral representation, one needs the $\L^1$ version, which we deduce using change of measure technique.

We will need the Burkholder-Davis-Gundy inequality (see \cite{Mey}) (for $p=1$) which states that there exist universal
 constants $c^{1}, c^{2}$ such that for all martingales $M$ and $T<\infty$,
\[
c^{1}\Eb  [ ([M,M ]_T)^{\frac{1}{2}} ]\le \Eb  [\sup_{0\le t\le T}  \lvert M_t \rvert ]\le c^{2}\Eb  [  ([M,M ]_T)^{\frac{1}{2}}].\]
After proving Theorem \ref{aztm1}, in the next section we will introduce sigma-martingales and give some elementary properties. Then we come to the main theorem on integral representation of martingales. This is followed by the second fundamental theorem of asset pricing.

\section{Proof of Theorem \ref{aztm1}}

We begin with a few auxiliary results.
In this section, we fix martingales $M^1,M^2,\ldots ,M^d$.
\bl\label{azl0} Let
{\small
\[\C_b(M^1,\ldots , M^d)=\{Z\in\M\,:\,Z_t=Z_0+\textstyle\sum_{j=1}^d\int_0^t f^j_s\d M^j_s,\;f^j\text{ bounded predictable },1\le j\le d\},\]
\[\F_b(M^1,\ldots , M^d)=\{Z\in\M\,: Z_t=Z_0+\textstyle\int_0^t f_s\d Y_s, \,f\in\L^1_m(Y),\,Y\in  \C_b(M^1,\ldots , M^d)\}.\]}
Then $\F_b(M^1,\ldots , M^d)=\F(M^1,\ldots , M^d)$.
\el
\pr Since bounded predictable process belong to $\L^1_m(N)$ for every martingale $N$, it follows that $\C_b(M^1,\ldots , M^d)\subseteq \C(M^1,\ldots , M^d)$ and hence $\F_b(M^1,\ldots , M^d)\subseteq \F(M^1,\ldots , M^d)$.

For the other part, let $Z\in \F$ be given by
\[Z_t= Z_0+\int_0^t f_s\d Y_s , \,f\in\L^1_m(Y)\]
where 
\[Y_t=\sum_{j=1}^d\int_0^t g^j_s\d M^j_s\]
with $g^j\in\L^1_m(M^j)$. Let $\xi_s=1+\sum_{j=1}^d\lvert g^j_s\rvert$, $h^j_s=\frac{g^j_s}{\xi_s}$ and
\[V_t=\sum_{j=1}^d\int_0^t h^j_s\d M^j_s.\]
Since $h^j$ are bounded, it follows that $V\in \C_b(M^1,M^2,\ldots , M^d)$.
Using $g^j_s=\xi_sh^j_s$ and $g^j\in\L^1_m(M^j)$, it follows that $\xi\in\L^1_m(V)$ and
\[Y_t=\int_0^t\xi_s\d V_s.\]
Since $f\in\L^1_m(Y)$, it follows that $f\xi\in\L^1_m(V)$ and
$\int f\d Y=\int f\xi \d V$. 
\qed

\bl \label{azl1}Let $Z\in\M$ be such that there exists a sequence of stopping times $\sigma_k\uparrow\infty$ with $\Eb_\Pb[\sqrt{[Z,Z]_{\sigma_k}}]<\infty$ and $X^k\in \F(M^1,M^2,\ldots , M^d)$ where $X^k_t=Z_{t\wedge\sigma_k}$. Then $Z\in\F(M^1,M^2,\ldots , M^d)$.
\el
\pr Let $X^k=\int f^k\d Y^k$ for $k\ge 1$ with $f^k\in\L^1_m(Y^k)$ and $Y^k\in\C_b(M^1,M^2,\ldots , M^d)$. Let $\phi^{k,j}$ be bounded predictable processes such that
\[Y^k_t=\sum_{j=1}^d\int_0^t\phi^{k,j}_s\d M^j_s.\]
Let $c_k>0$ be a common bound for $\phi^{k,1},\phi^{k,2},\ldots ,\phi^{k,d}$. Let us define $\eta^j,f$ by
\[\eta^j_t=\sum_{k=1}^\infty \frac{1}{c_k}\phi^{k,j}_t\Ind_{\{(\sigma_{k-1},\sigma_k]\}}(t).\]
\[f_t=\sum_{k=1}^\infty c_kf^k_t\Ind_{\{(\sigma_{-1},\sigma_k]\}}(t).\]
\[Y_t=\sum_{j=1}^d \int_0^t \eta^j_s\d M^j_s .\]
By definition, $\eta^j$ is bounded by 1 for every $j$ and thus
 $Y\in \C_b(M^1,M^2,\ldots , M^d)$.
We can note that
\[\begin{split}
Z_{t\wedge\sigma_k}- Z_{t\wedge\sigma_{k-1}}&=X^k_{t\wedge\sigma_k}-X^k_{t\wedge\sigma_{k-1}}\\
&=\int_0^tf^k_s\Ind_{\{(\sigma_{k-1},\sigma_k]\}}(s)\d Y^k_s\\
&=\int_0^tf_s\Ind_{\{(\sigma_{k-1},\sigma_k]\}}(s)\d Y_s.\end{split}
\]
Thus
\[Z_{t\wedge\sigma_k}=Z_0+\int_0^t\Ind_{[0,{\sigma_k}]}(s)f_s\d Y_s\]
and hence
\[[Z,Z]_{\sigma_k}=\int_0^{\sigma_k}(f_s)^2\d [Y,Y]_s.\]
Since by assumption, for all $k$
\[\Eb_\Pb[\sqrt{[Z,Z]_{\sigma_k}}\;]<\infty\]
it follows that $f\in\L^1_m(Y)$.
This proves the required result.
\qed
\bl \label{azl2}Let $Z^n\in\M$ be such that $\Eb[\lvert Z^n_t-Z_t\rvert ]\rightarrow 0$ for all $t$. Then there exists a sequence of stopping times $\sigma_k\uparrow\infty$ and a subsequence $\{n^j\}$ such that  each $k\ge 1$,
 \[\Eb[\sqrt{[Z,Z]_{\sigma_k }}]<\infty\] and writing $Y^j=Z^{n^j}$,
\be\label{az1}\Eb[\sqrt{[Y^j-Z,Y^j-Z]_{\sigma_k }}\;]\rightarrow 0 \text{ as }j\uparrow\infty.\ee
\el
\pr Let $n_0=0$. For each $j$, $\Eb[\lvert Z^n_j-Z_j\rvert ]\rightarrow 0$ as $n\rightarrow\infty$ and hence we can choose $n^j>n^{(j-1)}$ such that
\[\Eb[\lvert Z^{n^j}_j-Z_j\rvert ]\le 2^{-j}.\]
Then using Doob's maximal inequality we have
\[\Pb(\sup_{t\le j}\lvert Z^{n^j}_t-Z_t\rvert \ge \frac{1}{j^2})\le \frac{j^2}{2^j}.\]
As a consequence, writing $Y^j=Z^{n^j}$, we have
\be\label{az31}
\eta_t=\sum_{j=1}^\infty \sup_{s\le t} \lvert Y^j_s-Z_s\rvert<\infty \;\;a.s. \text{ for all }t<\infty.\ee
Now define
\[U_t=\lvert Z_t\rvert+\sum_{j=1}^\infty \lvert Y^j_t-Z_t\rvert.\]
In view of \eqref{az31}, it follows that $U$ is r.c.l.l. adapted process. For any stopping time $\tau\le m$, we have
\[\begin{split}
\Eb[U_\tau]&= \Eb[\lvert Z_\tau\rvert]+ \sum_{j=1}^\infty \Eb[\lvert Y^j_\tau-Z_\tau\rvert]\\
&\le \Eb[\lvert Z_m\rvert]+\sum_{j=1}^\infty \Eb[\lvert Y^j_m-Z_m\rvert]\\
&\le \Eb[\lvert Z_m\rvert]+\sum_{j=1}^m \Eb[\lvert Y^j_m-Z_m\rvert]+\sum_{j=m+1}^\infty 2^{-j} \\
&<\infty.
\end{split}\]
Here, we have used that $Z$, $Y^j-Z$ being martingales, $\lvert Z\rvert$, $\lvert Y^j-Z\rvert$ are sub-martingales and $\tau\le m$. Now defining 
\[\sigma_k=\inf\{t\,:\;U_t\ge k\text{ or } U_{t-}\ge k\}\wedge k\]
it follows that $\sigma_k$ are bounded stop times increasing to $\infty$ with
\[
\sup_{s\le \sigma_k}  U_s \le k+ U_{\sigma_k}
\]
and hence
\be\label{az31w}
\Eb[\sup_{s\le \sigma_k}  U_s ]<\infty.
\ee
Thus, for each $k$ fixed $\Eb[\sup_{s\le \sigma_k} \lvert Z_s\rvert]<\infty$ and thus 
Burkholder-Davis-Gundy inequality ($p=1$ case), we have $\Eb[\sqrt{[Z,Z]_{\sigma_k }}]<\infty$. 
In view of \eqref{az31} \[\lim_{j\rightarrow\infty}\sup_{s\le \sigma_k} \lvert Y^j_s-Z_s\rvert= 0\;\;\;\;a.s.\]
 and is dominated by $(\sup_{s\le \sigma_k}  U_s) $ which in turn is integrable as seen in \eqref{az31w}. Thus by dominated convergence theorem, we have
\[
\lim_{j\rightarrow\infty}\Eb[\sup_{s\le \sigma_k} \lvert Y^j_s-Z_s\rvert]= 0. \]
 The result \eqref{az1} now follows from Burkholder-Davis-Gundy inequality ($p=1$ case).
\qed

\bl \label{azl3}Let $V\in\F(M^1,M^2,\ldots , M^d)$ and $\tau$ be a bounded stopping time such that $\Eb[\sqrt{[V,V]_{\tau }}\;]<\infty$. Then for $\epsilon>0$, there exists $U\in 
\C_b(M^1,M^2,\ldots , M^d)$ such that
\[ \Eb[\sqrt{[V-U,V-U]_{\tau }}\;]\le\epsilon.\]
\el
\pr Let $V=\int f\d X$ where $f\in\L^1_m(X)$ and $X\in\C_b(M^1,M^2,\ldots , M^d)$. Since
\[[V,V]_t=\int_0^t\lvert f_s\rvert^2\d[X,X]_s,\]
the assumption on $V$ gives
\be\label{az33}
\Eb[\sqrt{\textstyle\int_0^\tau\lvert f_s\rvert^2\d[X,X]_s }]<\infty.
\ee
Defining $f^k_s=f_s\Ind_{\{\lvert f_s\rvert\le k\}}$, let
\[U^k=\int f^k\d X.\]
Since $X\in \C_b(M^1,M^2,\ldots , M^d)$ and $f^k$ is bounded, it follows that 
\[U^k\in\C_b(M^1,M^2,\ldots , M^d).\]
Note that as $k\rightarrow\infty$,
\[\Eb[\sqrt{[V-U^k,V-U^k]_{\tau }}]=\Eb[\sqrt{\textstyle\int_0^\tau\lvert f_s\rvert^2\Ind_{\{\lvert f_s\rvert> k\}}\d[X,X]_s }]\rightarrow 0\]
in view of \eqref{az33}. The result now follows by taking $U=U^k$ with $k$ large enough so that $\Eb[\sqrt{[V-U^k,V-U^k]_{\tau }}]<\epsilon.$

\qed

\bl \label{azl4} Suppose $Z\in \M$ and $\tau$  is a bounded stopping time such that $\Eb[\sqrt{[Z,Z]_{\tau }}\;]<\infty$, $Z_t=Z_{t\wedge\tau}$.
Let
$U^n\in \C_b(M^1,M^2,\ldots , M^d)$  with $U^n_0=0$ be such that
\[\Eb[\sqrt{[U^n-Z,U^n-Z]_{\tau }}\;]\le 4^{-n}.\]
Then there exists $X\in\C_b(M^1,M^2,\ldots , M^d)$ and $f\in\L^1_m(X)$ such that
\be\label{az66}Z_t=Z_0+\int_0^t f_s\d X_s.\ee 
\el
\pr 
Since $U^n\in\C_b(M^1,M^2,\ldots , M^d)$ with $U^n_0=0$, get bounded predictable processes $\{f^{n,j}:n\ge 1,\,1\le j\le d\}$ such that
\be\label{az2}U^n_t=\sum_{j=1}^d\int_0^tf^{n,j}_s\d M^j_s.\ee
Without loss of generality, we assume that $U^n_t=U^n_{t\wedge\tau}$ and $f^{n,j}_s=f^{n,j}_s\Ind_{[0,\tau]}(s)$.
Let 
\[\zeta=\sum_{n=1}^\infty 2^n\sqrt{[U^n-Z,U^n-Z]_{\tau }}.\]
Then $\Eb[\zeta]<\infty$ and hence $\Pb(\zeta<\infty)=1.$ Let
\[\eta=\zeta+\sqrt{[Z,Z]_{\tau }}+\sum_{j=1}^d\sqrt{[M^j,M^j]_{\tau }}\]
Let $c=\Eb[\exp\{-\eta\}]$ and let $\Qb$ be the probability measure on $(\Omega,\clf)$ defined by
\[\frac{\d\Qb}{\d\Pb}=\frac{1}{c}\exp\{-\eta\}.\]
Then it follows that $\alpha=\Eb_\Qb[\eta^2]<\infty$. Noting that
\[\eta^2\ge [Z,Z]_{\tau }+\sum_{j=1}^d[M^j,M^j]_{\tau }+\sum_{n=1}^\infty 2^{2n}[U^n-Z,U^n-Z]_{\tau }\]
we have $\Eb_\Qb[[Z,Z]_{\tau }]<\infty$, $\Eb_\Qb[[M^j,M^j]_{\tau }]<\infty$ for $1\le j\le d$. Likewise, $\Eb_\Qb[[U^n-Z,U^n-Z]_{\tau }]<\infty$ and so   $\Eb_\Qb[[U^n,U^n]_{\tau }]<\infty$. Note that $Z, M^j$ are no longer a martingales under $\Qb$, but we do not need that. 

Let $ \widetilde{\Omega}=[0,\infty)\times\Omega$. 
 Recall  that the predictable $\sigma$-field $\clp$ is the smallest $\sigma$ field on $ \widetilde{\Omega}$ with respect to which all continuous adapted processes are measurable. We will define signed measures $\Gamma_{ij}$ on  $\clp$ as follows: for $E\in\clp$, $1\le i,j\le d$ let
\[\Gamma_{ij}(E)=\int_{\Omega}\int_0^\tau \Ind_E(s,\omega)\d [M^i,M^j]_s(\omega)\d\Pb(\omega).\]
 Let $\Lambda=\sum_{j=1}^d\Gamma_{jj}$. From the properties of quadratic variation $[M^i,M^j]$, it follows that for all $E\in \clp$, the matrix $((\Gamma_{ij}(E)))$ is non-negative definite. Further, $\Gamma_{ij}$ is absolutely continuous w.r.t. $\Lambda$ $\forall i,j$. It follows (see appendix) that we can get predictable processes $c^{ij}$ such that
\[\frac{\d\Gamma_{ij}}{\d\Lambda}=c^{ij}\]
and that $C=((c^{ij}))$ is a non-negative definite matrix. By construction $\lvert c^{ij}\rvert \le 1$. We can diagonalise $C$ (i.e. obtain singular value decomposition) in a measurable way (see appendix) to obtain predictable processes $b^{ij}$, $d^j$ such that for all $i,k$, (writing $\delta_{ik}=1$ if $i=k$ and $\delta_{ik}=0$ if $i\neq k$),)
\be\label{az5}\sum_{j=1}^db^{ij}_sb^{kj}_s=\delta_{ik}\ee
\be\label{az6}\sum_{j=1}^db^{ji}_sb^{jk}_s=\delta_{ik}\ee
\be\label{az7}\sum_{j,l=1}^db^{ij}_sc^{jl}_sb^{kl}_s=\delta_{ik}d^{i}_s\ee 
%\be\label{az8}\sum_{j=1}^db^{ij}_sd^j_sb^{kj}_s=c^{ij}_s.\ee 
Since $((c^{ij}_s))$ is non-negative definite, it follows that $d^i_s\ge 0$. For $1\le j\le d$, 
let 
\[N^k=\sum_{l=1}^d\int b^{kl}_s\d M^l.\] Then $N^k$ are $\Pb$- martingales since $b^{ik}$ is a bounded predictable process. Indeed, $N^k\in\C_b(M^1,M^2,\ldots , M^d)$. Further, for $i\neq k$
\[[N^i,N^k]=\sum_{j,l=1}b^{ij}_sb^{kl}_s\d [M^j,M^l]\]
and hence for any bounded predictable process $h$
\[\begin{split}
\Eb_\Qb[\int_0^\tau h_s\d[N^i,N^k]]&=\int_\Omega\int_0^\tau h_s\sum_{j,l=1}^db^{ij}_sb^{kl}_s\d [M^j,M^l]\d\Pb(\omega)\\
&=\int_{\bar{\Omega}}h\sum_{j,l=1}^db^{ij}b^{kl}\d\Gamma_{jl}\\
&=\int_{\bar{\Omega}}h\sum_{j,l=1}^db^{ij}b^{kl}c^{jl}\d\Lambda\\
&=0\end{split}
\]  
where the last step follows from \eqref{az7}. As a consequence, for bounded predictable $h^i$,
\be\label{az10}
\Eb_\Qb[\sum_{i,k=1}^d\int_0^\tau h^i_sh^k_s\d[N^i,N^k]_s]=\Eb_\Qb[\sum_{k=1}^d\int_0^\tau (h^k_s)^2\d[N^k,N^k]_s]
\ee
Let us observe that \eqref{az10} holds for any predictable processes $\{h^i: 1\le i\le d\}$ provided the RHS if finite: we can first note that it holds for $\tilde{h}^i=h^i\Ind_{\{\lvert h\rvert\le c\}}$ where $\lvert h\rvert=\sum_{i=1}^d\lvert h^i\rvert$ and then let $c\uparrow\infty$. 
Note that for $n\ge m$
\[\begin{split}
\sqrt{[U^n-U^m,U^n-U^m]_{\tau }}&\le \sqrt{[U^n-Z,U^n-Z]_{\tau }}+ \sqrt{[U^m-Z,U^m-Z]_{\tau }}\\
&\le 2^{-m}\eta\end{split}
\]
and hence
\be\label{az11}
\Eb_\Qb[[U^n-U^m,U^n-U^m]_{\tau }]\le 4^{-m}\alpha.\ee
Let us define $g^{n,k}=\sum_{j=1}^df^{n,j}b^{kj}$. Then note that
\be\label{az12}
\begin{split}
\sum_{k=1}^d\int g^{n,k}\d N^k&=\sum_{k=1}^d\sum_{j=1}^d\int f^{n,j}b^{kj}\d N^k\\
&=\sum_{k=1}^d\sum_{j=1}^d\sum_{l=1}^d\int f^{n,j}b^{kj}b^{kl}_s\d M^l\\
&=\sum_{j=1}^d\int f^{n,j}\d M^j\\
&=U^n \end{split}\ee
where in the last but one step, we have used \eqref{az6}. Noting that  for $n\ge m$,
\be\label{az13}
\Eb_\Qb[\,[U^n-U^m,U^n-U^m]_\tau\,]=\Eb_\Qb[\sum_{k=1}^d\int_0^\tau(g^{n,k}_s-g^{m,k}_s)^2\d[N^k,N^k]_s]\ee
and using \eqref{az11}, we conclude
\be\label{az14}\begin{split}
\Qb(\sum_{k=1}^d\int_0^\tau(g^{n,k}_s-g^{m,k}_s)^2\d[N^k,N^k]\ge \frac{1}{m^4})&\le m^4\Eb_\Qb[\,[U^n-U^m,U^n-U^m]_\tau\,]\\
&\le \alpha m^44^{-m}.\end{split}\ee
Since $\Eb_\Qb[\,[M^i,M^i]_\tau]<\infty$ for all $i$ and $g^{n,i}$ are bounded, it follows that $\Eb_\Qb[\,[N^i,N^i]_\tau]<\infty$. Thus defining a measure  $\Theta$ on $\clp$ by
\[\Theta(E)=\int[\sum_{k=1}^d\int_0^\tau\Ind_E(s,\omega)\d[N^k,N^k]_s(\omega)\d\Qb(\omega)]\]
we get (using \eqref{az11} and \eqref{az13})
\[\int (g^{m+1,k}-g^{m,k})^2\d\Theta\le \alpha 4^{-m}\]
and as a consequence, using Cauchy-Schwartz inequality, we get
\[\int \sum_{m=1}^\infty\lvert g^{m+1,k}-g^{m,k}\rvert \d\Theta\le \sqrt{\Theta(\bar{\Omega})\alpha}<\infty.\]
Defining 
\[g^k_s=\limsup_{m\rightarrow\infty}g^{m,k}_s\]
it follows that $g^{m,k}\rightarrow g^k$ a.s. $\Theta$ and as a consequence, 
taking limit in \eqref{az14} as $n\rightarrow \infty$, we get
\be\label{az15}
\Qb(\sum_{k=1}^d\int_0^\tau(g^{k}_s-g^{m,k}_s)^2\d[N^k,N^k]_s\ge \frac{1}{m^4})
\le m^44^{-m}.\ee
Since $\Qb$ and $\Pb$ and  equivalent measures, it follows that 
\be\label{az16}
\Pb(\sum_{k=1}^d\int_0^\tau(g^{k}_s-g^{m,k}_s)^2\d[N^k,N^k]_s\ge \frac{1}{m^4})
\rightarrow 0 \text{ as }m\rightarrow\infty.\ee
In view of \eqref{az12}, we have for $m\le n$
\be\label{az3}[U^n,U^n]_\tau=\sum_{i,j=1}^d\int_0^\tau g^{n,i}_sg^{n,j}_s\d [N^i,N^j]_s\ee
and
\be\label{az4}[U^n-U^m,U^n-U^m]_\tau=\sum_{i,j=1}^d\int_0^\tau (g^{n,i}_s-g^{m,i}_s)(g^{n,j}_s-g^{m,j}_s)\d [N^i,N^j]_s.\ee
Taking limit in \eqref{az3} as $n\rightarrow\infty$, we get (using Fatou's lemma) 
\be\label{az20}
\Eb_\Pb[\sqrt{\textstyle\sum_{i,j=1}^d\int_0^\tau g^{i}_sg^{j}_s\d [N^i,N^j]}\,]\le \Eb_\Pb[\,\sqrt{[Z,Z]_\tau}\, ]\ee
(since \eqref{az1} implies $ \Eb_\Pb[\,\sqrt{[U^n,U^n]_\tau}\, ]\rightarrow \Eb_\Pb[\,\sqrt{[Z,Z]_\tau}\, ]$). 
Let us define bounded predictable processes $\phi^{j}$ and predictable process $h^n,h$ and a $\Pb$-martingale $X$ as follows:
\be\label{az21}
h_s=1+\sum_{i=1}^d\lvert g^{i}_s\rvert\ee
\be\label{az22}
\phi^{j}_s=\frac{g^{j}_s}{h_s}
\ee
\be\label{az23}
X_t=\sum_{j=1}^d\int_0^t\phi^{j}_s\d N^j_s\ee
Since $\phi^{j}$ is predictable, $\lvert \phi^{j}\rvert\le 1$ it follows that  $X\in\C_b(M^1,M^2,\ldots,M^d)$ and
\be\label{az24}
[X,X]_t=\sum_{j,k=1}^d\int_0^t\phi^{j}_s\phi^{k}_s\d [N^j,N^k]_s.\ee
Noting that $g^{j}_s=h_s\phi^{j}_s$ by definition, we conclude using \eqref{az20} that
\[
\int_0^t(h_s)^2\d [X,X]_s=\sum_{j,k=1}^d\int_0^tg^{j}_sg^{k}_s\d [N^j,N^k]_s
\]
and hence that 
\be\label{az25}
\Eb_\Pb[\sqrt{\textstyle\int_0^\tau(h_s)^2\d [X,X]_s}]\le \Eb_\Pb[\,\sqrt{[Z,Z]_\tau}\, ]\ee
Since $h=h\Ind_{[0,\tau]}$, we conclude that $h\in\L^1_m(X)$ and $Y=\int h\d X$ is a martingale with $Y_t=Y_{t\wedge\tau}$ for all $t$. Observe that
\[ [U^n,X]_t=\sum_{k,j=1}^d\int_0^tg^{n,k}_s\phi^j_s\d[N^k,N^j]_s\]
and hence
\[\begin{split}
[U^n,Y]_t&=\int_0^th_s\d[U^n,X]_s\\
&=\sum_{k,j=1}^d\int_0^tg^{n,k}_sh_s\phi^j_s\d[N^k,N^j]_s\\
&=\sum_{k,j=1}^d\int_0^tg^{n,k}_sg^j_s\d[N^k,N^j]_s\end{split}\]
and as a consequence
\[\begin{split}
[U^n-Y,U^n-Y]_t&=[U^n,U^n]_t-2[U^n,Y]_t+[Y,Y]_t\\
&=\sum_{k,j=1}^d\int_0^tg^{n,k}_sg^{n,j}_s\d[N^k,N^j]_s-2\sum_{k,j=1}^d\int_0^tg_s^{n,k}g^j_s\d[N^k,N^j]_s\\
&\;\;\;\;\;\;\;\;
+\sum_{k,j=1}^d\int_0^tg^{k}_sg^j_s\d[N^k,N^j]_s\\
&=\sum_{k,j=1}^d\int_0^t(g^{n,k}_s-g^k_s)(g^{n,j}_s-g^j_s)\d[N^k,N^j]_s\end{split}\]
and thus
\be\label{az26}
\Eb_\Qb[[U^n-Y,U^n-Y]_\tau]=\Eb_\Qb[\sum_{k=1}^d\int_0^\tau (g^{n,k}_s-g^k_s)^2\d[N^k,N^k]_s\ee
where we have used \eqref{az10}.
Taking $\liminf$ on the RHS in \eqref{az13} and using \eqref{az11}, we conclude
\[\Eb_\Qb[\sum_{k=1}^d\int_0^\tau (g^{k}_s-g^{m,k}_s)^2\d[N^k,N^k]_s]\le \alpha 4^{-m}\]
and hence
\[\Eb_\Qb[[U^n-Y,U^n-Y]_\tau]\le \alpha 4^{-n}.\]
Thus $[U^n-Y,U^n-Y]_\tau\rightarrow 0$ in $\Qb$-probability and hence in $\Pb$-probability. By assumption, $[U^n-Z,U^n-Z]_\tau\rightarrow 0$ in $\Pb$-probability.
Since
\[[Y-Z,Y-Z]_\tau\le 2([Y-U^n,Y-U^n]_\tau+[Z-U^n,Z-U^n]_\tau)\]
for every $n$, it follows that
\be\label{az27}
[Y-Z,Y-Z]_\tau=0\;\;a.s.\;\Pb.\ee
Since $Y,Z$ are $\Pb$-martingales such that $Z_t=Z_{t\wedge\tau}$ and $Y_t=Y_{t\wedge\tau}$,  \eqref{az27} implies $Y_t-Y_0=Z_t-Z_0$ for all $t$. 
Recall that by construction, $Y_0=0$, $Y=\int h\d X$ and $h\in\L^1_m(X)$, $X\in\C_b(M^1,M^2,\ldots , M^d)$. Thus \eqref{az66} holds.
\qed 

We now come to the proof of Theorem \ref{aztm2}. Let $Z^n\in \F(M^1,M^2,\ldots , M^d)$ be such that $\Eb[\lvert Z^n_t-Z_t\rvert ]\rightarrow 0$ for all $t$. We have to show that $Z\in \F(M^1,M^2,\ldots , M^d)$. 

Using Lemma \ref{azl2}, get a sequence of stopping times $\sigma_k\uparrow\infty$ and a subsequence $\{n^j\}$ such that $Y^j=Z^{n^j}$ satisfies for each $k\ge 1$, $\Eb[\sqrt{[Z,Z]_{\sigma_k }}]<\infty$ and 
\[\Eb[\sqrt{[Y^j-Z,Y^j-Z]_{\sigma_k }}\;]\rightarrow 0 \text{ as }j\uparrow\infty.\]

Let us now fix a $k$ and let $\tilde{Z}_t=Z_{t\wedge\sigma_k}$. We will show that $\tilde{Z}\in\F(M^1,M^2,\ldots , M^d)$. This will complete the proof in view of Lemma \ref{azl1}. 

Using Lemma \ref{azl3}, for each $n$, get $j_n$ such that  
\[\Eb[\sqrt{[Y^{j_n}-\tilde{Z},Y^{j_n}-\tilde{Z}]_{\sigma_k }}\;]\le 4^{-j-1}\]
For each $n$, taking $V=Y^{j_n}$ and $\epsilon=4^{-j-1}$, get $U^n\in \C_b(M^1,M^2,\ldots , M^d)$ such that
\[ \Eb[\sqrt{[Y^{j_n}-U^n,Y^{j_n}-U^n]_{\sigma_k }}\;]\le 4^{-j-1}.\]

Then we have
\[ \Eb[\sqrt{[U^n-\tilde{Z},U^n-\tilde{Z}]_{\sigma_k }}\;]\le 4^{-j}.\]
with $U^n\in \C_b(M^1,M^2,\ldots , M^d)$. 

This $\tilde{Z}\in\F(M^1,M^2,\ldots , M^d)$ in view of Lemma \ref{azl4}
\qed 

Now we turn to proof of Theorem \ref{aztm1}. Let $\xi^n\in\G_T$ be such that $\xi^n\rightarrow\xi$ in $\L^1(\Omega,\clf, \Pb)$. Let $\xi^n=X^n_T$ where $X^n\in\F(M^1,M^2,\ldots , M^d)$. Let us define $Z^n_t=X^n_{t\wedge T}$. Then $Z^n\in\F(M^1,M^2,\ldots , M^d)$ and the assumption on $\xi^n$ implies
\[Z^n_t\rightarrow Z_t=\Eb[\xi\mid\clf_t]\text{ in }\L^1(\Omega,\clf, \Pb)\;\;\forall t.\]
Thus Theorem \ref{aztm2} implies $Z\in \F(M^1,M^2,\ldots , M^d)$ and thus $\xi=Z_T$ belongs to $\G_T$.
\qed
\section{Sigma-martingales}
For a semimartingale $X$, let $L(X)$ denote the class of predictable process $f$ 
 such that $X$ admits a decomposition $X=N+A$ with $N$ being a local martingale, $A$ being a process with finite variation paths with $f\in\L^1_m(N)$ and 
\be\label{ax1}
\int_0^t\lvert f_s\rvert \d\lvert A\rvert_s<\infty \;\;a.s.\;\;\forall t<\infty.\ee
Then for $f\in\L(X)$, the stochastic integral $\int f\d X$ is defined as $\int f\d N+\int f\d A$. It can be shown that the definition does not depend upon the decomposition $X=N+A$. See \cite{J78}.

Let  $M$ be a martingale,  $f\in\L(M)$ and $Z=\int f\d M$. Then $Z$ is a local martingale if and only if $f\in\L^1_m(M)$. In answer to a question raised by P. A. Meyer, Chou \cite{Chou} introduced a class  $\Sigma_m$ of semimartingales consisting of $Z=\int f\d M$ for $f\in\L(M)$. Emery \cite{Em80} constructed example of $f,M$ such that $f\in\L(M)$ but $Z=\int f\d M$ is not a local martingale. Such processes occur naturally in mathematical finance and have been called sigma-martingales by Delbaen and Schachermayer\cite{DS98}.  

\definition{A semimartingale $X$ is said to be a sigma-martingale if there exists a $(0,\infty)$ valued predictable process $\phi$ such that $\phi\in\L(X)$ and 
\be\label{ax3}M_t=\int_0^t \phi_s\d X_s\ee
is a martingale.}    
Our first observation is:
\bl\label{axl5} Every local martingale $N$ is a sigma-martingale.
\el
\pr Let $\eta_n\uparrow\infty$ be a sequence of stopping times such that $N_{t\wedge\eta_n}$ is a martingale,
\[\sigma_n=\inf\{t\ge 0\,:\;\lvert N_t\rvert\ge n\text{ or }\lvert N_{t-}\rvert\ge n\}\wedge n\]
and $\tau_n=\sigma_n\wedge\eta_n$, then it follows that $N_{t\wedge\tau_n}$ is a uniformly integrable martingale and $a_n=\Eb[\,[N,N]_{\tau_n}]<\infty$. 

Define
\[h_s=\sum_{n=1}^\infty\frac{1}{2^n(1+a_n)}\Ind_{(\tau_{n-1},\tau_n]}.\]
Then $h$ being bounded belongs to $\L(N)$ and $M=\int h\d N$ is a local martingale with 
\be\label{ax4}\sup_{t<\infty}\Eb[\,[M,M]_t]<\infty.\ee
Thus $M$ is a uniformly integrable martingale. Since $h$ is $(0,\infty)$ valued by definition, it follows that $N$ is a sigma-martingale.
\qed
This leads to
\bl\label{axl4} A semimartingale $X$ is a sigma-martingale if and only if there exists a uniformly integrable martingale $M$ satisfying \eqref{ax4} and a predictable process $\psi\in\L(M)$ such that
\be\label{ax5}X_t=\int_0^t \psi_s\d M_s.\ee
\el
\pr
%Note that if $X,\phi,M$ are as in \eqref{ax3} with $\phi$ being $(0,\infty)$ valued, then $\psi_s=\frac{1}{\phi_s}$ belongs to $\L(M)$ and \eqref{ax5} holds. Conversely, 
Let $X$ be given by \eqref{ax5} with $M$ being a martingale satisfying \eqref{ax4} and $\psi\in\L(M)$, then defining
\[g_s=\frac{1}{(1+(\psi_s)^2)}, \;\;N_t=\int_0^tg_s\d X_s\]
it follows that $N=\int g\psi\d M$. Since $g\psi$ is bounded by 1 and $M$ satisfies \eqref{ax4}, it follows that $N$ is a martingale. Thus $X$ is a sigma-martingale.

Conversely, given a sigma-martingale $X$ and a $(0,\infty)$ valued predictable process  
$\phi$ such that $N=\int \phi\d X$ is a martingale, get $h$ as in Lemma \ref{axl5} and let $M=\int h\d N=\int h\phi \d X$. Then $M$ is a uniformly integrable martingale that satisfies \eqref{ax4} and $h\phi$ is a $(0,\infty)$ valued predictable process.
\qed

From the definition, it is not obvious that sum of sigma-martingales is also a sigma-martingale, but this is so as the next result shows.
\bt Let $X^1,X^2$ be sigma-martingales and $a_1,a_2$ be real numbers. Then $Y=a_1X^1+a_2X^2$ is also a sigma-martingale. 
\et
\pr Let $\phi^1,\phi^2$ be $(0,\infty)$ valued predictable processes such that
\[M^i_t=\int_0^t\phi^i_s\d X^i_s,\;\;i=1,2\]
are uniformly integrable martingales. Then, writing $\xi=\min(\phi^1,\phi^2)$ and $\eta^i_s=\frac{\xi_s}{\phi^i_s}$, it follows that  
\[N^i_t=\int_0^t \eta^i_s\d M^i_s=\int_0^t\xi_s\d X^i_s\]
are uniformly integrable martingales since $\eta^i$ is bounded by one. Clearly, $Y=a_1X^1+a_2X^2$ is a semimartingale and $\xi\in\L(X^i)$ for $i=1,2$ implies $\xi\in\L(Y)$ and
\[\int_0^t\xi_s\d Y_s=a_1N^1_s+a_2N^2_s\]
is a uniformly integrable martingale. Since $\xi$ is $(0,\infty)$ valued predictable process, it follows that $Y$ is a sigma-martingale.
\qed
The following result gives conditions under which a sigma-martingale is a local martingale.
\bl\label{axl1} Let $X$ be a sigma-martingale with $X_0=0$. Then $X$ is a local martingale if and only if there exists a sequence of stopping times $\tau_n\uparrow\infty$ such that
\be\label{ax7}
\Eb[\,\sqrt{[X,X]_{\tau_n}}\,]<\infty\;\;\forall n.\ee
\el
\pr
Let $X$ be a sigma-martingale and $\phi,\psi,M$ be such that \eqref{ax3}, \eqref{ax4} holds. Let  $\psi_s=\frac{1}{\phi_s}$ and as noted above, \eqref{ax5} holds. Then 
\[[X,X]_t=\int_0^t(\psi_s)^2\d [M,M]_s.
\]
Defining $\psi^k_s=\psi_s\Ind_{\{\lvert \psi_s\rvert\le k\}}$, it follows that 
\[X^k=\int_0^t\psi^k_s\d M_s\]
is a uniformly integrable martingale. Noting that for $k\ge 1$
\[[X-X^k,X-X^k]_t=\int_0^t(\psi_s)^2\Ind_{\{k<\lvert \psi_s\rvert\}}\d [M,M]_s\]
the assumption \eqref{ax7} implies that for each $n$ fixed,
\[\Eb[\,\sqrt{[X-X^k,X-X^k]_{\tau_n}}\,]\rightarrow 0\text{ as }k\rightarrow\infty.\]
The Burkholder-Davis-Gundy inequality ($p=1$) implies
that for each $n$ fixed,
\[\Eb[\,\sup_{0\le t\le\tau_n}\lvert X_t-X^k_t\rvert\,]\rightarrow 0\text{ as }k\rightarrow\infty.\] 
and as a consequence  $X^{[n]}_t=X_{t\wedge\tau_n}$ is a martingale for all $n$ and so $X$ is a local martingale. Conversely, if $X$ is a local martingale with $X_0=0$, and $\sigma_n$ are stop times increasing to $\infty$ such that $X_{t\wedge\sigma_n}$ are martingales then defining $\zeta_n=\inf\{t\,:\,\lvert X_t\rvert\ge n\}$ and
$\tau_n=\sigma_n\wedge\zeta_n$, it follows that $\E[\lvert X_{\tau_n}\rvert]<\infty$ and since
\[\sup_{t\le \tau_n}\lvert X_t\rvert\le n+\lvert X_{\tau_n}\rvert\]
it follows that $\E[\sup_{t\le \tau_n}\lvert X_t\rvert]<\infty$. Thus, \eqref{ax7} holds in view of Burkholder-Davis-Gundy inequality ($p=1$).
\qed

The previous result gives us:
\bc \label{azc1} A bounded sigma-martingale $X$ is a martingale.
\ec
\pr Since $X$ is bounded, say by $K$, it follows that jumps of $X$ are bounded by $2K$. Thus jumps of the increasing process $[X,X]$ are bounded by $4K^2$ and thus $X$ satisfies \eqref{ax7} for
\[\tau_n=\inf\{t\ge 0\,:\;[X,X]_t\ge n\}.\]
Thus $X$ is a local martingale and being bounded, it is a martingale.
\qed
Here is a variant of the example given by Emery \cite{Em80} of a sigma-martingale that is not a local martingale.. Let $\tau$ be a random variable with exponential distribution (assumed to be $(0,\infty)$ valued without loss of generailty) and $\xi$ with $\Pb(\xi=1)=\Pb(\xi=-1)=0.5$, independent of $\tau$. Let 
\[M_t=\xi\Ind_{[\tau,\infty)}(t)\]
and $\clf_t=\sigma(M_s:s\le t)$. Easy to see that $M$ is a martingale. Let $f_t=\frac{1}{t}\Ind_{(0,\infty)}(t)$ and $X_t=\int_0^tf\d M$. Then  $X$ is a sigma-martingale and 
\[[X,X]_t=\frac{1}{\tau^2}\Ind_{[\tau,\infty)}(t).\]
For any stopping time $\sigma$, it can be checked that  $\sigma$ is a constant on $\sigma<\tau$ and thus if $\sigma$ is not identically equal to 0, $\sigma\ge (\tau\wedge a)$ for some $a>0$. Thus, $\sqrt{[X,X]_\sigma}\ge\frac{1}{\tau}\Ind_{\{\tau<a\}}$. It follows that for any stop time $\sigma$, not identically zero, $\Eb[\sqrt{[X,X]_\sigma}]=\infty$ and so $X$ is not a local martingale.

The next result shows that $\int f\d X$ is a sigma-martingale if $X$ is one.
\bl\label{axl7} Let $X$ be a sigma-martingale, $f\in\L(X)$ and let $U=\int f\d X$. Then $U$ is a sigma-martingale. 
\el
\pr Let $M$ be a martingale and $\psi\in\L(M)$ be such that $X=\int \psi \d M$ (as in Lemma \ref{axl4}). Now $U=\int f\d X=\int f\psi \d M$. Thus, once again invoking Lemma \ref{axl4}, one concludes that $X$ is a sigma-martingale.
\qed

We now introduce the class of equivalent sigma-martingale measures (ESMM) and show that it is a convex set.
\definition{Let $X^1,\ldots ,X^d$ be r.c.l.l. adapted processes and let $\E^s(X^1,\ldots ,X^d)$ denote the class of probability measures $\Qb$  such that $X^1,\ldots ,X^d$ are sigma-martingales w.r.t. $\Qb$.}
Let 
\[\E^s_\Pb(X^1,\ldots ,X^d)=\{\Qb\in \E^s(X^1,\ldots ,X^d)\,:\, \Qb\text{ is equivalent to }\Pb\}\]
and
\[\tilde{\E}^s_\Pb(X^1,\ldots ,X^d)=\{\Qb\in \E^s(X^1,\ldots ,X^d)\,:\, \Qb\text{ is absolutely continuous w.r.t. }\Pb\}\]

\bt \label{esmm} For semimartingales $X^1,\ldots ,X^d$, $\E^s(X^1,\ldots ,X^d)$, $\E^s_\Pb(X^1,\ldots ,X^d)$ and $\tilde{\E}^s_\Pb(X^1,\ldots ,X^d)$ are convex sets.
\et
\pr Let us consider the case $d=1$. Let $\Qb^1,\Qb^2\in\E^s(X)$. Let $\phi^1,\phi^2$ be $(0,\infty)$ valued predictable processes such that
\[M^i_t=\int_0^t \phi^i_s\d X_s\]
are martingales under $\Qb_i$, $i=1,2$. Let $\phi_s=min(\phi^1_s,\phi^2_s)$ and let 
\[M_t=\int_0^t\phi_s\d X_s.\]
Noting that $M_t=\int_0^t\xi_s^i\d M^i_s $ where $\xi_s={\phi_s}(\phi_s^i)^{-1}$ is bounded, it follows that $M$ is a martingale under $\Qb^i,i=1,2$. Now if $\Qb$ is any convex combination of $\Qb^1,\Qb^2$, it follows that $M$ is a $\Qb$ martingale and hence $X_t=\int_0^t(\phi_s)^{-1}\d M_s$ is a sigma-martingale under $\Qb$. Thus 
$\E^s_\Pb(X)$ is a convex set. Since
$\E^s(X^1,\ldots ,X^d)=\cap_{j=1}^d\E^s(X^j)$
it follows that $\E^s(X^1,\ldots ,X^d)$ is convex. Convexity of $\E^s_\Pb(X^1,\ldots ,X^d)$ and $\tilde{\E}^s_\Pb(X^1,\ldots ,X^d)$ follows from this.
\qed
In analogy with the definition of $\C$ for martingales $M^1,\ldots ,M^d$, 
for sigma-martingales $M^1,M^2,\ldots , M^d$  let
\[\C(M^1,M^2,\ldots , M^d)=\{Z\in\M\,:\,Z_t=Z_0+\sum_{j=1}^d\int_0^t f^j_s\d M^j_s,\;f^j\in\L^1_m(M^j)\}\]
\[\F(M^1,\ldots , M^d)=\{Z\in\M\,: Z_t=Z_0+\int_0^t f_s\d Y_s, \,f\in\L^1_m(Y),\,Y\in  \C(M^1,\ldots , M^d)\}.\]
\bl\label{ayl9} Let $M^1,\ldots, M^d$ be
sigma-martingales and let $\phi^j$ be $(0,\infty)$ valued predictable processes such that 
\be\label{ay9}N^j_t=\int_0^t\phi^j_s\d M^j_s\ee
are uniformly integrable martingales. Then
\be\label{ay10}
\C(M^1,M^2,\ldots , M^d)=\C(N^1,N^2,\ldots ,N^d)
\ee
\be\label{ay10a}
\F(M^1,M^2,\ldots , M^d)=\F(N^1,N^2,\ldots ,N^d).
\ee
\el
\pr
Let $\psi^j_s=(\phi^j_s)^{-1}$. Note that $M^j=\int \psi^j\d N^j$. 
Then if $Y\in \C(M^1,M^2,\ldots , M^d)$ is given by
\be\label{ay11}
Y_t=\sum_{j=1}^d\int_0^tf^j_s\d M^j_s,\;\;f^j\in\L(M^j)\ee
then defining $g^j=f^j\psi^j$, we can see that $g^j\in\L(N^j)$ and $\int f^j\d M^j=\int g^j\d N^j$. Thus
\be\label{ay12}
Y_t=\sum_{j=1}^d\int_0^tg^j_s\d N^j_s,\;\;g^j\in\L(N^j).\ee
Similarly, if $Y\in\C(N^1,N^2,\ldots ,N^d)$ is given by \eqref{ay12}, then defining $f^j=\phi^jg^j$, we can see that $Y$ satisfies \eqref{ay11}. Thus \eqref{ay10} is true. Now \eqref{ay10a} follows from \eqref{ay10}.
\qed
\section{Integral Representation w.r.t. Martingales}
Let $M^1,\ldots, M^d$ be sigma-martingales. 
\definition{ A sigma-martingale $N$ is said to have an integral representation w.r.t. 
$M^1,\ldots, M^d$ if $N\in \F(M^1,M^2,\ldots , M^d)$ or in other words,
$\exists Y\in  \C(M^1,M^2,\ldots , M^d)$ and  $f\in\L(Y)$ such that
\be\label{ay21}
N_t=N_0+\int_0^tf_s\d Y_s\,\;\;\forall t.
\ee
}

Here is another observation needed later.
\bl\label{cal1} 
Let $M$ be a $\Pb$-martingale. Let $\Qb$ be a probability measure equivalent to $\Pb$. Let $\xi=\frac{\d\Pb}{\d\Qb}$ and let $Z$ be the r.c.l.l. martingale given by $Z_t=\Eb_\Pb[\xi\mid\clf_t]$. Then
\begin{enumerate}[(i)]
\item $M$ is a $\Qb$-martingale if and only if $MZ$ is a $\Pb$-martingale.
\item $M$ is a $\Qb$-local martingale if and only if $MZ$ is a $\Pb$-local martingale.
\item If $M$ is a $\Qb$-local martingale then $[M,Z]$ is a $\Pb$-local martingale.
\item If $M$ is a $\Qb$-sigma-martingale then $[M,Z]$ is a $\Pb$-sigma-martingale.
\end{enumerate} 
\el
\pr For a stopping time $\sigma$, 
let $\eta$ be a non-negative $\clf_\sigma$ measurable random variable. Then 
\[\Eb_\Qb[\eta]=\Eb_\Pb[\eta Z]=\Eb_\Pb[\eta \Eb[Z\mid\clf_\sigma]\,]=\Eb_\Pb[\eta Z_\sigma].\]
Thus $M_s$ is $\Qb$ integrable if and only if $M_sZ_s$ is $\Pb$-integrable.
 Further, for 
any stopping time $\sigma$,
\be\label{bm60}
\Eb_\Qb[M_\sigma]=\Eb_\Pb[M_\sigma Z_\sigma].\ee
Thus (i) follows from the observation that an integrable adapted process $N$ is a martingale if and only if $\Eb[N_\sigma]=\Eb[N_0]$ for all bounded stopping times $\sigma$.  For (ii), if $M$ is a $\Qb$-local martingale, then get stopping times $\tau_n\uparrow\infty$ such that for each $n$, $M_{t\wedge\tau_n}$ is a martingale. Then we have
\be\label{bm61}\Eb_\Qb[M_{\sigma\wedge\tau_n}]=
\Eb_\Pb[M_{\sigma\wedge\tau_n} Z_{\sigma\wedge\tau_n}].
\ee
Thus $M_{t\wedge\tau_n}Z_{t\wedge\tau_n}$ is a $\Pb$-martingale and thus $MZ$ 
is a $\Pb$- local martingale. The converse follows similarly.

For (iii), note that $M_tZ_t=M_0Z_0+\int_0^tM_{s-}\d Z_s+\int_0^tZ_{s-}\d M_s +[M,Z]_t$ and the two stochastic integrals are $\Pb$ local martingales, the result follows from (ii). For (iv), representing the $\Qb$ sigma-martingale $M$ as $M=\int \psi \d N$, where $N$ is a $\Qb$ martingale and $\psi\in\L(N)$, we see
\[[M,Z]=\int_0^t\psi_s\d [N,Z]_s.\]
By (iii), $[N,Z]$ is a $\Qb$ sigma-martingale and hence $[M,Z]$ is a  $\Qb$ sigma-martingale.
\qed

The main result on integral representation is:
\bt \label{intrep} Let $\clf_0$ be trivial. Let $M^1,\ldots, M^d$ be sigma-martingales on $(\Omega,\clf,\Pb)$. Then the following are equivalent.
\begin{enumerate}[(i)]
\item Every bounded martingale admits representation w.r.t. $M^1,\ldots, M^d$.
\item Every uniformly integrable martingale admits representation w.r.t. $M^1,\ldots, M^d$.
\item Every sigma-martingale admits representation w.r.t. $M^1,\ldots, M^d$.
\item $\Pb$ is an extreme point of the convex set $\E^s(M^1,\ldots, M^d)$. 
\item $\tilde{\E}^s_{\Pb}(M^1,\ldots, M^d)=\{\Pb\}.$
\item $\E^s_{\Pb}(M^1,\ldots, M^d)=\{\Pb\}.$
\end{enumerate}
\et
\pr  Since every bounded martingale is uniformly integrable and a uniformly integrable martinagle is a sigma-martingale, we have\\
\centerline{ (iii)$\Rightarrow$ (ii) $\Rightarrow$ (i).}

\noindent (i) $\Rightarrow$ (ii) is an easy consequence of Theorem \ref{aztm2}: given a uniformly integrable martingale $Z$, for $n\ge 1$, let us define martingales $Z^n$ by
\[Z^n_t=\Eb[Z\Ind_{\{\lvert Z\rvert\le n \}}]\mid \clf_t].\]
We take the r.c.l.l. version of the martingale. It is easy to see that $Z^n$ are bounded martingales and in view of (i), $Z^n\in\F(M^1,M^2,\ldots , M^d)$. Moreover, for $n\ge t$
\[\Eb[\lvert Z^n_t-Z_t\rvert ]\le \Eb[Z\Ind_{\{\lvert Z\rvert> n \}}]\]
and hence for all $t$, $\Eb[\lvert Z^n_t-Z_t\rvert ]\rightarrow 0$. Theorem \ref{aztm2} now implies $Z\in\F(M^1,M^2,\ldots , M^d)$. This proves (ii).

\noindent We next prove (ii) $\Rightarrow$ (iii). Let $X$ be a sigma-martingale. In view of Lemma \ref{axl4}, get a uniformly integrable martingale $N$ and a predictable process $\psi$ such that 
\[X=\int \psi \d N.\]
Let $N_t=N_0+\int_0^t f_s\d Y_s$ where $Y\in\C(M^1,M^2,\ldots , M^d)$. Then  we have
\[X_t=X_0+\int_0^t f_s\psi\d Y_s\]
and thus $X$ admits an integral representation w.r.t. $M^1,\ldots, M^d$.

Suppose (v) holds and suppose  $\Qb_1,\Qb_2$ $\in\E^s(M^1,M^2,\ldots , M^d)$ and $\Pb=\alpha \Qb_1+(1-\alpha)\Qb_2$. It follows that $\Qb_1,\Qb_2$ are absolutely continuous w.r.t. $\Pb$ and hence $\Qb_1, \Qb_2\in \tilde{\E}_{\Pb}^s(M^1,M^2,\ldots , M^d)$. In view of (v), $\Qb_1=\Qb_2=\Pb$ and thus $\Pb$ is an extreme point of $\E^s(M^1,M^2,\ldots , M^d)$ and so (iv) is true. 

Since $\E_{\Pb}^s(M^1,M^2,\ldots , M^d) \subseteq \tilde{\E}_{\Pb}^s(M^1,M^2,\ldots , M^d)$, it follows that (v) implies (vi). On the other hand, suppose (vi) is true and $\Qb\in\tilde{\E}_{\Pb}^s(M^1,M^2,\ldots , M^d)$. Then $\Qb_1=\frac{1}{2}(\Qb+\Pb)\in \E_{\Pb}^s(M^1,M^2,\ldots , M^d)$. Then (vi) implies $\Qb_1=\Pb$ and hence $\Qb=\Pb$ and thus (v) holds.

 Till now we have proved (i) $\Longleftrightarrow$ (ii) $\Longleftrightarrow$ (iii) and  
(iv) $\Leftarrow$ (v) $\Longleftrightarrow$ (vi). To complete the proof, we will show (iii) $\Rightarrow$ (v) and (iv) $\Rightarrow$ (i).

Suppose that (iii) is true and let $\Qb\in\tilde{\E}^s_\Pb(M^1,M^2,\ldots , M^d)$. Now let $\xi$ be the Radon-Nikodym derivative of $\Qb$ w.r.t. $\Pb$. 
Let $R$ denote the r.c.l.l. martingale: $R_t=\Eb[\xi\mid\clf_{t}]$.  Since $\clf_0$ is trivial, $N_0=1$. In view of property (iii),   we can get  $Y\in\C(M^1,M^2,\ldots , M^d)$
and a predictable processes $f\in\L(Y)$ such that 
\be\label{ax51}
R_t=1+\int_0^tf_s\d Y_s.
\ee
Note that
\be\label{ax51d} [R,R]_t=\int_0^tf_s^2\d [Y,Y]_s.\ee
Since $M^j$ is a sigma-martingale under $\Qb$ for each $j$, it follows that $Y$ is a $\Qb$ sigma-martingale. By Lemma \ref{cal1}, this gives $[Y,R]$ is a $\Pb$ sigma-martingale and hence
\be\label{ax51a}
V^k_t=\int_0^tf_s\Ind_{\{\lvert f_s\rvert\le k\}}\d [Y,R]_s\ee
is a $\Pb$ sigma-martingale. Noting that
\[ [Y,R]_t=\int_0^tf_s\d [Y,Y]_s\]
we see that
\be\label{ax51b}
V^k_t=\int_0^tf^2_s\Ind_{\{\lvert f_s\rvert\le k\}}\d [Y,Y]_s.\ee
Thus we can get $(0,\infty)$ valued predictable processes $\phi^j$ such that
\[U^k_t=\int_0^t\phi^k_s\d V^k_s\]
is a martingale. But $U^k$ is a non-negative martingale with $U^k_0=0$. As a result $U^k$ is identically equal to 0 and thus so is $V^k$. It then follows that (see \eqref{ax51d}) $[R,R]=0$ which yields $R$ is identical to 1 and so $\Qb=\Pb$. Thus   
$\tilde{\E}^s_\Pb(M^1,M^2,\ldots , M^d)$ is a singleton.
Thus (iii) $\Rightarrow$ (v). 

To complete the proof, we will now prove that (iv) $\Rightarrow$ (i).  
Suppose $M^1,M^2,\ldots , M^d$ are such that $\Pb$ is an extreme point of $\E^s(M^1,M^2,\ldots , M^d)$. Since $M^j$ is a sigma-martingale under $\Pb$, we can choose $(0,\infty)$ valued predictable $\phi^j$ such that
\[N^j_t=\int_0^t\phi^j_s\d M^j_s\]
is a uniformly integrable martingale under $\Pb$ and as seen in Lemma \ref{ayl9}, we then have
\[\F(M^1,M^2,\ldots , M^d)=\F(N^1,N^2,\ldots ,N^d).\]
Suppose (i) is not true. We will show that this leads to a contradiction.
So suppose $S$ is a bounded martingale that does not admit representation w.r.t. $
M^1,M^2,\ldots , M^d$, {\em i.e.} $S\not\in \F(M^1,M^2,\ldots , M^d)=\F(N^1,N^2,\ldots ,N^d)$, then for some $T$, 
\[S_T\not\in \K_T(N^1,N^2,\ldots ,N^d)\]
We have proven in Theorem \ref{aztm1} that $ \K_T(N^1,N^2,\ldots ,N^d)$ is closed in $\L^1(\Omega,\clf, \Pb)$.
Since $\K_T$ is not equal to $\L^1(\Omega, \clf_T,\Pb)$, by the Hahn-Banach Theorem, there exists $\xi\in\L^\infty(\Omega, \clf_T,\Pb)$, $\Pb(\xi\neq 0)>0$  such that 
\[\int \eta\xi \d\Pb=0\;\;\forall \eta\in\K_T.\]
Then for all constants $c$, we have
\be\label{azk44}
\int \eta(1+c\xi) \d\Pb=\int\eta\d\Pb\;\;\forall \eta\in\K_T.\ee
Since $\xi$ is bounded, we can choose a $c>0$ such that
\[\Pb(c\lvert\xi\rvert<\frac{1}{2})=1.\]
Now, let $\Qb$ be the measure with density $\eta=(1+c\xi)$. Then $\Qb$ is a probability measure.  Thus \eqref{azk44} yields
\be\label{azk45}
\int \eta\d\Qb=\int\eta\d\Pb\;\;\forall \eta\in\K_T.\ee
For any bounded stop time $\tau$and $1\le j\le d$, $N^j_{\tau\wedge T}\in\K_T$ and hence
\be\label{azk46}
\Eb_\Qb[N^j_{\tau\wedge T}]=\Eb_\Pb[N^j_{\tau\wedge T}]=N^j_0\ee
On the other hand, 
\be\label{azk47}\begin{split}
\Eb_\Qb[N^j_{\tau\vee T}]&=\Eb_\Pb[\eta N^j_{\tau\vee T}]\\
&=\Eb_\Pb[\Eb_\Pb[\eta N^j_{\tau\vee T}\mid \clf_T]]\\
&=\Eb_\Pb[\eta\Eb_\Pb[ N^j_{\tau\vee T}\mid \clf_T]]\\
&=\Eb_\Pb[\eta N_T]\\
&=\Eb_\Qb[ N_T]\\
&=N^j_0.\end{split}\ee
where we have used the facts that $\eta$ is $\clf_T$ measurable, $N^j$ is a $\Pb$ martingale and \eqref{azk46}. Now
\[ \Eb_\Qb[N^j_{\tau}]=\Eb_\Qb[N^j_{\tau\wedge T}]+\Eb_\Qb[N^j_{\tau\vee T}]-\Eb_\Qb[N^j_{T}]=N^j_0.\]
Thus $N^j$ is a $\Qb$ martingale and since
\[M^j=\int_0^t\frac{1}{\phi^j_s}\d N^j_s\]
it follows that $M^j$ is a $\Qb$ sigma-martingale. 
Thus $\Qb\in  \E^s(M^1,\ldots, M^d)$.
Similarly, if $\tilde{\Qb}$ is the measure with density $\eta=(1-c\xi)$, we can prove that $\tilde{\Qb}\in  \E^s(M^1,\ldots, M^d)$. Since $\Pb=\frac{1}{2}(\Qb+\tilde{\Qb})$, this contradicts the assumption that $\Pb$ is an extreme point of  $\E^s(M^1,\ldots, M^d)$. 
Thus (iv) $\Rightarrow$ (i). This completes the proof.
\qed 
\section{Completeness of Markets}
Let the (discounted) prices of $d$ securities be given by $X^1,\ldots ,X^d$. We assume that $X^j$ are semimartingales and that they satisfy the property NFLVR so that an ESMM exists. 
\bt \label{sftap} {\bf The Second Fundamental Theorem Of Asset Pricing} \\Let $X^1,\ldots ,X^d$ be semimartingales on $(\Omega,\clf,\Pb)$ such that $\E^s_\Pb(X^1,\ldots, X^d)$ is non-empty. Then the following are equivalent:
\begin{enumerate}[(a)]
\item For all $T<\infty$, for all $\clf_T$ measurable bounded random variables $\xi$ (bounded by say $K$), there exist $g^j\in\L(X^j)$ with
\be\label{aza1}
Y_t=\sum_{j=1}^d\int_0^tg^j_s\d X^j_s\ee
a constant $c$ and  $f\in\L(Y)$ such that $\lvert \int_0^tf_s\d Y_s\rvert \le 2K$  and 
\be\label{aza2}
\xi=c+\int_0^Tf_s\d Y_s.\ee
 \item The set $\E^s_\Pb(X^1,\ldots, X^d)$ is a singleton. 
\end{enumerate}
\et
\pr First suppose that $\E^s_\Pb(X^1,\ldots, X^d)=\{\Qb\}.$ Consider the martingale
$M_t=\Eb_\Qb[\xi\mid\clf_t]$. Note that $M$ is bounded by $K$. In view of the equivalence of (i) and (v) in Theorem \ref{intrep}, we get that $M$ admits a representation w.r.t. $X^1,\ldots, X^d$ - thus we get  $g^j\in\L(X^j)$ and $f\in\L(Y)$ where $Y$ is given by by \eqref{aza1}, with
\[M_t=M_0+\int_0^tf_s\d Y_s.\]
Since $\clf_0$ is trivial, $M_0$ is a constant. Since $M$ is bounded by $K$, it follows that $\int_0^tf_s\d Y_s$ is bounded by $2K$. Thus (b) implies (a).
 
Now suppose (a) is true. Let $\Qb$ be an ESMM. Let $M_t$ be a martingale. We will show that $M\in\F(X^1,\ldots, X^d)$, {\em i.e.} $M$ admits integral representation w.r.t. $X^1,\ldots, X^d$. In view of Lemmas \ref{azl1} and \ref{ayl9}, suffices to show that for each $T<\infty$, $N\in\F(X^1,\ldots, X^d)$, where $N$ is  defined by $N_t=M_{t\wedge T}$. 

Let $\xi=N_T$. Then in view of assumption (a), we have
\[\xi=c+\int_0^Tf_s\d Y_s\]
with $Y$ given by \eqref{aza1}, a constant $c$ and  $f\in\L(Y)$ such that $U_t=\int_0^tf_s\d Y_s$ is  bounded. Since $U$ is a sigma-martingale that is bounded, it follows that $U$ is a martingale. It follows that
\[N_t=c+\int_0^tf_s\d Y_s, \;\;0\le t\le T.\]
Thus $N\in \F(X^1,\ldots, X^d)$. 

We have proved that (i) in Theorem  \ref{intrep} holds and hence (v) holds, {\em i.e.} the ESMM is unique.
\qed

\appendix
\begin{center}
      {\bf APPENDIX}
    \end{center}
  \renewcommand{\theequation}{A.\arabic{equation}}
  % redefine the command that creates the equation no.
  \setcounter{equation}{0}  % reset counter 
 % \section*{APPENDIX}  % use *-form to suppress numbering
   
For a non-negative definite symmetric matrix $C$, the eigenvalue-eigenvector decomposition gives us a representation
\be\label{apx1}C=B^TDB\ee
\be\label{apx2}\mbox{$B$ is a orthogonal matrix and $D$ is a diagonal matrix.}\ee
This decomposition is not unique, but for each non-negative definite symmetric matrix $C$, the set of pairs $(B,D)$ satisfying \eqref{apx1}-\eqref{apx2} is compact. Thus it admits a measurable selection, in other words, there exists a Borel mapping $\theta$ such that $\theta(C)=(B,D)$ where $B,C,D$ satisfy \eqref{apx1}-\eqref{apx2}.
(See \cite{Graf} or Corollary 5.2.6 \cite{SMS}).

Let $\cld$ be a $\sigma$-field on a non-empty set $\Gamma$  and for $1\le i,j\le d$, $\lambda_{ij}$ be $\sigma$-finite signed measures on $(\Gamma, \cld)$ such that\\ 
\centerline{ For all $E\in\cld$, the matrix$((\lambda_{ij}(E)))$ is a symmetric non-negative definite matrix.}
Let $\theta(E)=\sum_{i=1}^d\lambda_{ii}(E)$. Then  for $1\le i,j\le d$ there exists a version $c^{ij}$ of the Radon-Nikodym derivate $\frac{d\lambda_{ij}}{d\theta}$ such that for all $\gamma\in\Gamma$, the matrix $((c^{ij}(\gamma)))$ is non-negative definite.

To see this, for $1\le i\le j\le d$ let $f^{ij}$ be a version of the Radon-Nikodym derivative $\frac{d\lambda_{ij}}{d\theta}$ and let $f^{ji}=f^{ij}$. For rationals $r_1,r_2,\ldots , r_d$, let 
\[A_{r_1,r_2,\ldots , r_d}=\{\gamma:\sum_{ij}r_ir_jf^{ij}(\gamma)< 0\}.\]
 Then $\theta(A_{r_1,r_2,\ldots , r_d})=0$ and hence $\theta(A)=0$ where
 \[A=\cup\{A_{r_1,r_2,\ldots , r_d}:\;   r_1,r_2,\ldots , r_d\text{ rationals}\}.\]
 The required version is now given by
\[c^{ij}(\gamma)=f^{ij}(\gamma)\Ind_{A^c}(\gamma).\]

\vfill

\noindent {\it Address}:  H1 Sipcot IT Park, Siruseri, Kelambakkam 603103, India\\
{\it E-mail}: rlk@cmi.ac.in, bvrao@cmi.ac.in
\end{document}